\pdfoutput=1
\RequirePackage{ifpdf}
\ifpdf 
\documentclass[pdftex]{sigma}
\else
\documentclass{sigma}
\fi

\numberwithin{equation}{section}

\newtheorem{Theorem}{Theorem}[section]
\newtheorem{Lemma}[Theorem]{Lemma}
\newtheorem{Proposition}[Theorem]{Proposition}
 { \theoremstyle{definition}
\newtheorem{Definition}[Theorem]{Definition}
\newtheorem{Remark}[Theorem]{Remark} }

\usepackage{tikz}

\usepackage{enumerate}

\begin{document}

\allowdisplaybreaks

\renewcommand{\thefootnote}{$\star$}

\newcommand{\arXivNumber}{1503.00169}

\renewcommand{\PaperNumber}{072}

\FirstPageHeading

\ShortArticleName{(Co)isotropic Pairs in Poisson and Presymplectic Vector Spaces}

\ArticleName{(Co)isotropic Pairs in Poisson\\ and Presymplectic Vector Spaces\footnote{This paper is a~contribution to the Special Issue
on Poisson Geometry in Mathematics and Physics.
The full collection is available at \href{http://www.emis.de/journals/SIGMA/Poisson2014.html}{http://www.emis.de/journals/SIGMA/Poisson2014.html}}}

\Author{Jonathan LORAND~$^\dag$ and Alan WEINSTEIN~$^\ddag$}

\AuthorNameForHeading{J.~Lorand and A.~Weinstein}

\Address{$^\dag$~Department of Mathematics, ETH Zurich, Zurich, Switzerland}
\EmailD{\href{mailto:jonathan.lorand@math.uzh.ch}{jonathan.lorand@math.uzh.ch}}
\URLaddressD{\url{http://www.math.uzh.ch/?assistenten&key1=11723}}

\Address{$^\ddag$~Department of Mathematics, University of California, Berkeley, CA 94720 USA}
\EmailD{\href{mailto:alanw@math.berkeley.edu}{alanw@math.berkeley.edu}}
\URLaddressD{\url{http://math.berkeley.edu/~alanw/}}

\ArticleDates{Received March 01, 2015, in f\/inal form September 03, 2015; Published online September 10, 2015}

\Abstract{We give two equivalent sets of invariants which classify pairs of coisotropic subspaces of  f\/inite-dimensional Poisson vector spaces. For this it is convenient to dualize; we work with pairs of isotropic subspaces of presymplectic vector spaces. We identify ten elementary types which are the building blocks of such pairs, and we write down a matrix, invertible over~$\mathbb{Z}$, which takes one 10-tuple of invariants to the other. }

\Keywords{coisotropic subspace; direct sum decomposition; Poisson vector space; presymplectic vector space}

\Classification{15A21; 18B10; 53D99}

\renewcommand{\thefootnote}{\arabic{footnote}}
\setcounter{footnote}{0}

\section{Introduction}

The classif\/ication of pairs of linear coisotropic subspaces arises in many contexts, including the following two.

The f\/irst is the problem of classifying, up to conjugation by linear symplectomorphisms, canonical relations (lagrangian correspondences) from a f\/inite-dimensional symplectic vector space to itself.   Without symplectic structure, this
classif\/ication of linear relations was carried out by Towber~\cite{Towber} and is a special case of results of Gel'fand and Ponomarev~\cite{Gelfand_Ponomarev}.  For graphs of symplectomorphisms, the classif\/ication amounts to identifying the conjugacy classes in the group of symplectic matrices. This classif\/ication and the problem of f\/inding associated normal forms has a long history extending from Williamson~\cite{Williamson} to Gutt~\cite{Gutt}. In the general symplectic case, a result of Benenti and Tulczyjew~\cite[Proposizioni~4.4 \&~4.5]{Benenti_Tulczyjew} tells us that any canonical relation $X\leftarrow Y$ is given by coisotropic subspaces of $X$ and $Y$ and a symplectomorphism between the corresponding reduced spaces.  When $X = Y$, one step in the classif\/ication up to conjugacy of canonical relations is to classify the coisotropic pairs consisting of  the range and domain. The further steps of the classif\/ication remain as work in progress.   Partial results may be found in the f\/irst author's
Master's thesis~\cite{Lorand}.

The second context is that of extending the Wehrheim--Woodward theory of linear canonical relations (see \cite{Li-Bland_Weinstein,Weinstein}) to the case where the category of lagrangian correspondences between symplectic vector spaces is enlarged to the category of coisotropic correspondences.  Each pair of coisotropic subspaces of $X$ gives a WW morphism represented by a~diagram having the form
${\mathbf 1} \leftarrow X \leftarrow {\mathbf 1}$, and isomorphic pairs represent the same WW morphism. There are also nonisomorphic pairs representing the same WW morphism. The issue is to determine exactly which pairs are ``WW equivalent''.   This problem is now solved, as part of a complete description of the WW categories of (co)isotropic relations (see \cite{Weinstein2}).

So far, we have been discussing symplectic ambient spaces. But since co\-isotropic correspondences are fundamental in Poisson geometry\footnote{See, e.g., \cite{we:coisotropic}, where these correspondences are called Poisson relations.},
it is natural to consider the endomorphism classif\/ication and WW problems for linear coisotropic correspondences between general Poisson vector spaces.   This leads immediately to the classif\/ication problem for coisotropic pairs in Poisson vector spaces, which is the subject of this paper\footnote{The isomorphism classif\/ication of Poisson relations is treated in
\cite{lo-we:decomposition}.}.

It turns out to be simpler to work indirectly with Poisson vector spaces  via their duals, which are presymplectic (i.e., equipped with   possibly degenerate skew-symmetric bilinear forms).  The coisotropic subspaces are then replaced by their annihilators, which are isotropic.   This duality gives an equivalence between the Poisson/coisotropic and presymplectic/isotropic categories. In the following, therefore, we will concentrate on pairs of isotropic subspaces in presymplectic vector spaces. To classify such pairs, we f\/irst show in Section~\ref{decomposition}  that any isotropic pair can be decomposed as the direct sum of pairs, each of which is of one of ten ``elementary'' types.  In Section~\ref{indecomposables} we show that these elementary types are in turn isomorphic to multiples of ten indecomposable pairs
in spaces of dimension~1,~2, and~3, and that the multiplicities are invariants which fully classify isotropic pairs, up to conjugation by linear presymplectomorphisms. To prove that the multiplicities are invariant, we assemble them into a 10-vector and show that this vector is related by an invertible integer matrix to a 10-vector, each of whose entries is the dimension of  a space constructed in a simple way from a presymplectic space and an ordered pair of its isotropic subspaces.

The preliminary version~\cite{Lorand_Weinstein} of this paper contains our results in the symplectic case. We have in the meantime found the reference~\cite{Sergeichuk}, which treats subspace pairs in spaces with bilinear forms, using methods and results related to quiver representations. The present paper gives an elementary and essentially self-contained treatment  for spaces carrying a skew-symmetric form which is possibly degenerate; it remains for the future to relate our results to a broader representation-theoretic framework (see Remark~\ref{rep theory}).

\section{Preliminaries}\label{preliminaries}

Throughout this paper, $V$ will denote a f\/inite-dimensional presymplectic vector space, with presymplectic structure~$\omega$. For any subspace $W \subseteq V$, we call the subspace $\{ v \in V \,|\, \omega (v, w) = 0$ $\forall \,w \in W \}$ the \textit{orthogonal of~$W$} and denote it by~$W^{\perp}$.
 For the radical $V^{\perp}$ of $V$ we reserve the letter $R$. A~\textit{presymplectomorphism} from $V$ to~$\hat V$ is
is a linear isomorphism $\varphi\colon   V \rightarrow \hat{V}$  which pulls back the presymplectic structure on~$\hat V$ to the one on~$V$.

A subspace $A$ of $V$ is \textit{isotropic} if $A \subseteq A^{\perp}$.  An \textit{isotropic pair} in~$V$ is an ordered pair of isotropic subspaces in~$V$. Isotropic pairs $(A,B)$ and $(\hat A, \hat B)$  in $V$ and $\hat V$ respectively are \textit{equivalent} if there exists a presymplectomorphism $\varphi\colon V \rightarrow \hat V$ such that $\varphi(A) = \hat A$ and $\varphi(B) = \hat B$. In the Poisson setting, where a coisotropic subspace is a subspace annihilated by an isotropic in the dual, this equivalence corresponds to there being an invertible Poisson map which takes one coisotropic pair to the other. In the symplectic situation, when $\omega$ is non-degenerate, any coisotropic subspace is the orthogonal of an isotropic subspace. Clearly, a linear symplectomorphism will take one coisotropic pair to the other if and only if it maps the corresponding isotropic orthogonals to one another.

An isotropic pair $(A,B)$ in $V$ is the \textit{direct sum} of isotropic pairs $(A_1,B_1)$ and $(A_2,B_2)$ in~$V_1$ and~$V_2$, respectively, if $V = V_1 \oplus V_2$, $A = A_1 \oplus A_2$, $B = B_1 \oplus B_2$, and~$V_1$ is orthogonal to $V_2$ with respect to $\omega$.\footnote{In other words, $\omega$ is the direct sum of the presymplectic structures on~$V_1$ and~$V_2$.} In this case, we have $(A_i,B_i) = (A \cap V_i, B \cap V_i)$, for $i = 1,2$. The def\/inition of direct sum of isotropic pairs naturally extends to any f\/inite number of summands.

A decomposition
\begin{gather}\label{direct sum}
V = \bigoplus_{i=1}^m V_i,
\end{gather}
is \textit{orthogonal} if its summands are pairwise orthogonal; it is \textit{distributive}\footnote{It was pointed out to us by one of the referees that the existing language of graded vector spaces could be used here, but we f\/ind the ``distributivity'' terminology more convenient, because we make no use of any algebraic structure on the index set, and because it is superf\/luous to index the components in some decompositions. Of course, the basic properties associated to our def\/inition correspond to basic facts about graded vector spaces.}
 with respect to a~subspace $W \subseteq V$ if
\begin{gather*}
W = \bigoplus_{i=1}^m W \cap V_i.
\end{gather*}
The main task in Section~\ref{decomposition} will be to construct an orthogonal decomposition of $V$ which is distributive with respect to isotropic subspaces $A$ and $B$, and such that  each of the resulting summands has a simple form.

\begin{Remark}\label{modular iteration} A general strategy for constructing a direct sum decomposition of $V$ which is distributive with respect to $W \subseteq V$ is the following stepwise procedure.
In each step, f\/ind a~subspace $V' \subseteq V$ such that one of the following holds:
\begin{enumerate}[(i)]\itemsep=0pt
\item $W \cap V' = 0$,
\item $W \cap V' = V'$.
\end{enumerate}
If (i) is the case, then there exists a subspace $C \subseteq V$ such that $V = V' \oplus C$ and $W \subseteq C$. This decomposition of $V$ is, by construction, distributive with respect to $W$. If (ii) is the case, i.e., $V' \subseteq W$, then for any complement $C$ of $V'$ in $V$ one has $W = V' \oplus (W \cap C)$ by the modular law\footnote{The modular law is the fact that, for subspaces $E$, $F$ and $G$ of a vector space $V$, if $E \subseteq G$, then  $G \cap (E + F) = E + (G \cap F)$. See \cite[p.~56]{Roman}, for example.}, and so the resulting decomposition $V = V' \oplus C$ is distributive with respect to $W$.
\end{Remark}

The following lemma shows that this procedure will indeed achieve the desired result. If in each step $V'$ and $C$ can be chosen to be orthogonal, then the resulting decomposition will also be orthogonal.

\begin{Lemma}\label{distributive filter}
Suppose that  $V$ has a decomposition \eqref{direct sum}, and let $W \subseteq V$ be a subspace.
For each $l \in \{0,1,\dots ,m-1\}$ set
\begin{gather*}
C_l := \bigoplus_{i=l+1}^m V_i.
\end{gather*}
\begin{enumerate}[{\rm (i)}]\itemsep=0pt
\item If, for each $l \in \{1,\dots ,m-1\}$, the decomposition $V_l \oplus C_{l}$ is distributive with respect to $W \cap C_{l-1}$, then the decomposition~\eqref{direct sum} is distributive with respect to $W$.
\item If, for each $l \in \{1,\dots ,m-1\}$, the decomposition $V_l \oplus C_{l}$ is orthogonal, then the decomposition~\eqref{direct sum} is orthogonal.
\end{enumerate}
\end{Lemma}

\begin{proof}
(i) We apply the assumptions successively to construct a decomposition of $W$ composed of the intersections of $W$ with the $V_i$. For $l=1$, by assumption we have a decomposition
\begin{gather*}
W = W \cap V_1 \oplus W \cap C_1.
\end{gather*}
The assumption for $l=2$, applied to the second summand of this decomposition, gives
\begin{gather*}
W = W \cap V_1 \oplus W \cap V_2 \oplus W \cap C_2.
\end{gather*}
Clearly, proceeding in this manner for increasing $l$ will, after $m-1$ steps, lead to a decomposition
\begin{gather*}
W = W \cap V_1 \oplus W \cap V_2 \oplus  \dots  \oplus W \cap V_{m-1} \oplus W \cap C_{m-1},
\end{gather*}
which, after substitution using the identity $C_{m-1} = V_m$, is the desired result.

(ii) Choose any two indices $i, j \in \{1,\dots ,m\}$ such that $i \neq j$. We need to show that~$V_i$ and~$V_j$ are orthogonal. Because this relation is symmetric with respect to~$i$ and~$j$, we may assume without loss of generality that $i < j$. Then $i \leq m-1$ and $V_j \subseteq C_i$. By assumption~$V_i$ is orthogonal to~$C_i$, so in particular~$V_i$ is orthogonal to~$V_j$.
\end{proof}

In the next lemma, we collect some basic properties of distributive decompositions. The proofs are, as above, straightforward; we leave them to the reader.

\begin{Lemma}\label{distributivity properties}
Let $V = \bigoplus_{i=1}^m V_i$ be a decomposition which is distributive with respect to subspaces $E \subseteq V$ and $F \subseteq V$. Then,
\begin{enumerate}[{\rm (i)}]\itemsep=0pt
\item The decomposition is distributive with respect to $E \cap F$ and $E + F$, and
\begin{gather*}
(E + F) \cap V_i = (E \cap V_i) + (F \cap V_i)
\end{gather*}
for each~$i$.
\item If the decomposition of $V$ is orthogonal, then it is also distributive with respect to $E^{\perp},$ and
\begin{gather*}
E^{\perp} \cap V_i = (E \cap V_i)^{\perp} \cap V_i
\end{gather*}
for each $i$.
\end{enumerate}
\end{Lemma}

\begin{Remark}\label{radical distributivity}
As a special case of Lemma~\ref{distributivity properties}~(ii), with $V$ in the role of $E$, it follows that any orthogonal decomposition of $V$ is distributive with respect to the radical~$R$.
\end{Remark}

We also recall some basic facts from presymplectic linear algebra.

\begin{Lemma}\label{linear algebra}
Let $E$ and $W$ be subspaces of $V$.
\begin{enumerate}[{\rm (i)}]\itemsep=0pt
\item If $E$ is symplectic, i.e., $\omega \vert_{E}$ is non-degenerate, then $V = E \oplus E^{\perp}$.
\item If $E$ is such that $W = W \cap W^{\perp} \oplus E$, then $E$ is symplectic. In particular, any complement of~$R$ in~$V$ is symplectic.
\item $\dim W + \dim W^{\perp} =  \dim V + \dim W \cap R$.
\item $W^{\perp \perp} = W + R$.
\item When $V$ is symplectic, $W$ is  {lagrangian} in $V$, i.e., $W^{\perp} = W$, if and only if~$W$ is isotropic and $\dim W = \tfrac{1}{2} \dim V$.
\end{enumerate}
\end{Lemma}

We omit the proofs of these facts but use them to prove the following lemma which we will use later.

\begin{Lemma}\label{reduction}
Let $I \subseteq V$ be an isotropic subspace, and let $L$ be a complement of $I \cap R$ in $I$. If~$E$ is a complement of $R + I$ in $I^{\perp}$ and $L'$ is a~complement of $R + I$ in $E^{\perp}$, then $E$ and $L \oplus L'$ are symplectic subspaces, $L$ is a~lagrangian subspace of the latter, and $V = R \oplus E \oplus (L \oplus L')$.
\end{Lemma}

\begin{proof}
All references in this proof are to Lemma~\ref{linear algebra} above. From~(iv) one has that $I^{\perp \perp} = I + R$, and since $I$ is isotropic, $I + R \subseteq I^{\perp}$ (note that $R \subseteq I^{\perp}$ in any case). Thus, by~(ii), with $I^{\perp}$ in the role of~$W$, $E$ is symplectic. From this and (i) it follows that $V = E \oplus E^{\perp}$, and since by assumption $E^{\perp} = R \oplus L \oplus L',$ one obtains $V = R \oplus E \oplus L \oplus L'$. That $(L \oplus L')$ is symplectic follows from~(ii), with $E^{\perp}$ in the role of~$W$, since $E^{\perp} \cap (E^{\perp})^{\perp} = E^{\perp} \cap (I + R) = R$ and $E^{\perp} = R \oplus (L \oplus L')$. Finally, to show that $L$ is lagrangian in $(L \oplus L')$, one can apply~(iii), with~$L$ in the role of~$W$, which (using $L^{\perp} = I^{\perp} = R \oplus L \oplus E$ and $L \cap R = 0$) shows that $2 \dim L = \dim (L \oplus L')$. By~(v), this completes the proof.
\end{proof}

\section{Decomposition of isotropic pairs}\label{decomposition}

\begin{Proposition}\label{elementary decomposition}
Any isotropic pair $(A,B)$ in a presymplectic space $V$ can be decomposed as an orthogonal direct sum of ten isotropic pairs, each of which is of an \textit{elementary type}. If $(A_i, B_i)$ denotes the $i$-th summand, with ambient space~$V_i$, the elementary types are:
\begin{enumerate}\itemsep=0pt
\item[{\rm 1)}] $R_1 = V_1$, $A_1 = B_1=0$,
\item[{\rm 2)}] $R_2 = V_2$, $A_2 = B_2 =V_2$,
\item[{\rm 3)}] $R_3 = V_3$, $A_3 = V_3$ and $B_3=0$,
\item[{\rm 4)}] $R_4 = V_4$, $A_4 = 0$ and $B_4=V_4$,
\item[{\rm 5)}] $\dim V_5 = 3 \dim R_5$; $R_5$, $A_5$ and $B_5$ have pairwise zero intersection, and $A_5 \oplus R_5 = B_5 \oplus R_5 = A_5 \oplus B_5$,
\item[{\rm 6)}] $R_6 = 0$, $A_6$ and $B_6$ are lagrangian in $V_6$, and $A_6 = B_6$,
\item[{\rm 7)}] $R_7 = 0$, $A_7$ is lagrangian in $V_7$ and $B_7 = 0$,
\item[{\rm 8)}] $R_8 = 0$, $A_8 = 0$ and $B_8$ is lagrangian in $V_8$,
\item[{\rm 9)}] $R_9 = 0$, $A_9$ and $B_9$ are lagrangian in $V_9$, and $A_9 \cap B_9 = 0$,
\item[{\rm 10)}] $R_{10} = 0$, $A_{10} = B_{10} = 0$.
\end{enumerate}

We will refer to these types as type~$1$, type~$2$, etc.
\end{Proposition}

\begin{proof}
We will construct a ten-part orthogonal direct sum decomposition of $V$ by proceeding according to Remark~\ref{modular iteration}, successively peeling away subspaces in such a way that, in each step, the resulting decomposition is both orthogonal and distributive with respect to $A$ and $B$,  i.e., so that Lemma~\ref{distributive filter} will apply. At each step, we will f\/ind a summand $V_i$ such that $A_i = A \cap V_i$ and $B_i=B \cap V_i$ form a pair of type $i$ in $V_i$.

To begin, let $V_1$ be a complement of $R \cap (A + B)$ in $R$. Since $V_1 \cap (A + B) = 0$, we can choose a complement $C_1$ of $V_1$ in $V$ which contains $A + B$.   Since $V_1\subseteq R$, the decomposition $V=V_1 \oplus C_1$ is orthogonal.

Next, set $V_2 = A \cap B \cap R \cap C_1$ and let $C_2$ be any complement of $V_2$ in $C_1$ (in particular then $A = V_2 \oplus (A \cap C_2)$ and $B = V_2 \oplus (B \cap C_2)$).

Now, let $V_3 = A \cap R \cap C_2$; since $B \cap C_2 \cap V_3 = 0$, we can choose a complement $C_3$ of $V_3$ in $C_2$ such that $B \cap C_2 \subseteq C_3$. Set $V_4 = B \cap R \cap C_3$. We have $A = V_2 \oplus V_3 \oplus A \cap C_3$ and $A \cap V_4 = 0$; in particular we can choose a complement $C_4$ of $V_4$ in $C_3$ which contains $A \cap C_3 \subseteq C_4$.   Since $V_4 \subseteq B \cap C_3$, $B = V_2 \oplus V_4 \oplus B \cap C_4$. The decomposition $V = V_1 \oplus V_2 \oplus V_3 \oplus V_4 \oplus C_4$ is orthogonal, since the f\/irst four summands are subspaces of $R$.

$C_4$ is now small enough so that $R \cap C_4$ has zero intersection with each of the isotropics $A \cap C_4$ and $B \cap C_4$; we consider $Q_R = R \cap C_4 \cap (A + B)$, which is the radical of the presymplectic space~$C_4$. One might try to f\/ind a complement to $Q_R$ which contains both $A \cap C_4$ and $B \cap C_4$, but this is not in general possible. Instead, we consider the space $S = A \cap C_4 \cap (Q_R + B) + B \cap C_4 \cap (Q_R + A)$, which contains $Q_R$ and $A \cap B \cap C_4$ as independent subspaces. Let $Q$ be a complement of $A \cap B \cap C_4$ in $S$ such that $Q_R \subset Q$, and set $Q_A = A \cap C_4 \cap Q$, $Q_B = B \cap C_4 \cap Q$. Note that each of the summands in the def\/inition of $S$ decomposes as the direct sum of its intersections with $A \cap B \cap C_4$ and $Q$, and the latter subspaces span $Q$ and lie in $Q_A +Q_B$. This implies that $Q_A + Q_B = Q$. In fact, $Q = Q_A \oplus Q_B = Q_R \oplus Q_A  = Q_R \oplus Q_B$. To see this, it suf\/f\/ices to show that $Q_B \subseteq Q_A + Q_R$, since we know that $Q_R$, $Q_A$ and $Q_B$ are pairwise independent, so dimension considerations then give the result. If $v = r + a$ is an element of $Q_B \subseteq B \cap C_4 \cap (Q_R + A)$, with $r \in Q_R$ and $a \in A \cap C_4$, then $a = v - r$ must also lie in $Q$, since $v \in Q_B \subseteq Q$ and $r \in Q_R \subseteq Q$, which means that $a \in Q_A$, and thus $v \in Q_R \oplus Q_A$.

We will now use $Q$ to carry out the next step in our orthogonal decomposition of $V$.
First, note that $Q^{\perp} \cap C_4 = Q_A^{\perp} \cap C_4 = Q_B^{\perp} \cap C_4$; in particular, it follows that $Q^{\perp} \cap C_4$ contains $A \cap C_4 + B \cap C_4$. In this latter space, let $T$ be a complement of $A \cap B \cap C_4$ which contains $Q = Q_A \oplus Q_B$. Then, in particular, $A \cap C_4 = (A \cap B \cap C_4) \oplus (A \cap C_4 \cap T)$, with $Q_A \subseteq A \cap C_4 \cap T$, and $B \cap C_4 = (A \cap B \cap C_4) \oplus (B \cap C_4 \cap T)$, with $Q_B \subseteq B \cap C_4 \cap T$. Let $A'$ and $B'$ be complements of $Q_A$ and $Q_B$ in $A \cap C_4 \cap T$ and $B \cap C_4 \cap T$, respectively. Then, $A'$ and $B'$ are independent, and their direct sum is a complement of $Q$ in $T$. This gives a decomposition
\begin{gather}\label{sum decomp}
A \cap C_4 + B \cap C_4 = Q_A \oplus Q_B \oplus (A \cap B \cap C_4) \oplus A' \oplus B'.
\end{gather}
Now choose a complement $C_5$ of $Q$ in $Q^{\perp} \cap C_4$ which contains the last three summands in (\ref{sum decomp}), and let $P$ be a complement of $Q$ in $C_5^{\perp} \cap C_4$. By Lemma~\ref{reduction}, with $Q_A$ in the role of $I$ (recall that $Q_A^{\perp} \cap C_4 = Q^{\perp} \cap C_4)$, we have an orthogonal decomposition $C_4 = Q_R \oplus (Q_A \oplus P) \oplus C_5$, the spaces~$C_5$ and~$(Q_A \oplus P)$ are symplectic, and~$Q_A$ is lagrangian in the latter. We set $V_5 = Q_R \oplus Q_A \oplus P$. By construction, the decomposition $C_4 = V_5 \oplus C_5$ is distributive with respect to both $A \cap C_4$ and $B \cap C_4$: indeed $Q_A \subseteq A \cap V_5$ and $(A \cap B \cap C_4) \oplus A' \subseteq A \cap C_5$ form a decomposition of $A \cap C_4$, so the inclusions are in fact equalities; this is similarly the case for~$Q_B \subseteq B \cap V_5$ and $(A \cap B \cap C_4) \oplus B' \subseteq B \cap C_5$ in $B \cap C_4$.

At this point we have completely decomposed the radical $R$; it remains to decompose the symplectic space $C_5$ with respect to the isotropics $A \cap C_5$ and $B \cap C_5$. Let $C_6$ be a complement of $A \cap B \cap C_5$ in $(A \cap B)^{\perp} \cap C_5$; since $A \cap C_5$ and $B \cap C_5$ each contain $A \cap B \cap C_5$ and lie in $(A \cap B)^{\perp} \cap C_5$, they each can be decomposed as the sum of $A \cap B \cap C_5$ and their respective intersections with $C_6$.  Choose a complement of $A \cap B \cap C_5$ in $C_6^{\perp} \cap C_5$; by Lemma~\ref{reduction} it is non-degenerately paired with $A \cap B \cap C_5$, and the direct sum of these two spaces forms a symplectic space $V_6$ in which $A \cap B \cap C_5$ is lagrangian.

Next, consider $A \cap B^{\perp} \cap C_6$ and let $C_7$ be a complement of this subspace in $(A \cap B^{\perp})^{\perp} \cap C_6$. Since $A \cap B^{\perp} \cap B \cap C_6 = 0$ and $B \cap C_6 \subseteq (A \cap B^{\perp})^{\perp} \cap C_6$, we can choose $C_7$ so as to contain $B \cap C_6$.  By Lemma~\ref{reduction} again, any complement to $A \cap B^{\perp} \cap C_6$ in $C_7^{\perp} \cap C_6$ is non-degenerately paired with $A \cap B^{\perp} \cap C_6$, forming a symplectic subspace $V_7$ in which $A \cap B^{\perp} \cap C_6$ is a lagrangian subspace.

This same construction can now be applied in $C_7$, with the roles of $A$ and $B$ exchanged, to yield an orthogonal decomposition $C_7 = V_8 \oplus C_8$ into symplectic subspaces which is distributive with respect to $A \cap C_7$ and $B \cap C_7$ and is such that $B \cap A^{\perp} \cap C_7$ is a lagrangian subspace of~$V_8$.

In the symplectic space $C_8$, $A \cap C_8$ and $B \cap C_8$ are independent isotropic subspaces such that $A \cap B^{\perp} \cap C_8 = B \cap A^{\perp} \cap C_8 = 0$. We claim that $V_9 = (A \cap C_8) \oplus (B \cap C_8)$ and $V_{10} = A^{\perp} \cap B^{\perp} \cap C_8$ form an orthogonal decomposition $V_8 = V_9 \oplus V_{10}$ into symplectic subspaces. For this, it is suf\/f\/icient to show that $V_9$ and $V_{10}$ are independent, since $V_{10} = V_9^{\perp} \cap C_8$. So suppose $v \in V_9 \cap V_{10}$. As an element of $V_9$, $v = a + b$, with $a \in A \cap C_8$ and $ b \in B \cap C_8$. Since $A \cap C_8 \subseteq A^{\perp} \cap C_8$, $a$ is an element of $A^{\perp} \cap C_8$, as is~$v$. Thus $b = v - a$ lies in $B \cap A^{\perp} \cap C_8 = 0$. The same argument shows that $a =0$, and hence $v = 0$, as desired. Finally, we check that $A \cap C_8$ and $B \cap C_8$ are each lagrangian subspaces of $V_9$: on the one hand, being isotropic, they each may have at most half the dimension of~$V_9$; on the other hand, their direct sum spans~$V_9$, so they must each have dimension $\tfrac{1}{2}\dim V_9$,  i.e., they are each lagrangian.

This completes our proof, since by Lemma~\ref{distributive filter} the constructed decomposition of $V$ is both orthogonal and distributive with respect to~$A$ and $B$ (and also $R$,  cf.\ Remark~\ref{radical distributivity}), and it is clear from the construction that the intersections of these spaces with the summands $V_i$ each def\/ine isotropic pairs of the types stated in the proposition.
\end{proof}

\section{Indecomposables}\label{indecomposables}

We denote the canonical basis vectors in $\mathbb{R}^n$ by $e_1, e_2, \dots , e_n$. Angled brackets `{$\langle \  \rangle$}' indicate the span of a vector or a~list of vectors. The standard (pre)symplectic structure on $\mathbb{R}^n$ of rank $2k$ is the bilinear form given in canonical coordinates by the $n \times n$ matrix
\begin{gather*}
\left(
\begin{matrix}
0 & \mathbb{I}_k & 0  \\
-\mathbb{I}_k & 0 & 0  \\
0 & 0 & 0
\end{matrix}
\right),
\end{gather*}
where $\mathbb{I}_k$ denotes the $k \times k$ identity matrix, and the zeros denote square zero matrices of appropriate sizes. A~symplectic basis of a~vector space with symplectic form $\omega$ is one such that the coordinate matrix of $\omega$ with respect to this basis has the same form as the matrix of the standard symplectic form.

\begin{Definition}\label{indecomposable}
An isotropic pair $(A,B)$ in $V$ is \textit{indecomposable} if for any direct sum decomposition into isotropic pairs $(A_1,B_1)$ and $(A_2,B_2)$, corresponding to an orthogonal decomposition $V=V_1 \oplus V_2$, it follows that either $V_1 = 0$ or $V_2 = 0$.
\end{Definition}

\begin{Theorem}\label{classification}
Any indecomposable isotropic pair $(A,B)$ in a presymplectic space $V$ is isomorphic to one of the following isotropic pairs.

Pairs where the ambient space is $\mathbb{R}$, carrying the zero presymplectic structure:
\begin{enumerate}\itemsep=0pt
\item[{\rm 1)}]  $A = 0$, $B = 0$,
\item[{\rm 2)}]  $A = \langle e_1 \rangle$, $B = \langle e_1 \rangle$,
\item[{\rm 3)}]  $A = \langle e_1 \rangle$, $B = 0$,
\item[{\rm 4)}] $A = 0$, $B = \langle e_1 \rangle$.
\end{enumerate}

Pairs where the ambient space is $\mathbb{R}^3$, equipped with the standard presymplectic structure of rank $2$:
\begin{enumerate}\itemsep=0pt
\item[{\rm 5)}] $A = \langle e_1 \rangle$, $B = \langle e_1 + e_3 \rangle$.
\end{enumerate}

Pairs where the ambient space is $\mathbb{R}^2$, carrying the standard symplectic structure:
\begin{enumerate}\itemsep=0pt
\item[{\rm 6)}] $A = \langle e_1 \rangle$, $B = \langle e_1 \rangle$,
\item[{\rm 7)}] $A = \langle e_1 \rangle$, $B = 0$,
\item[{\rm 8)}] $A = 0$, $B = \langle e_2 \rangle$,
\item[{\rm 9)}] $A = \langle e_1 \rangle$, $B = \langle e_2 \rangle$,
\item[{\rm 10)}] $A = 0$, $B = 0$.
\end{enumerate}

Any isotropic pair can be decomposed as the direct sum of indecomposable isotropic pairs. The multiplicities $n_1, \dots  , n_{10}$ of the indecomposable types in the decomposition are invariants of the isotropic pair itself.  Two isotropic pairs are equivalent if and only if their corresponding multiplicities are equal.
\end{Theorem}

\begin{proof}
First of all, we check that each of the isotropic pairs listed above is indeed indecomposable. For types~1 through~4 and~6 through~10,  this follows directly from dimension considerations and, in the latter cases, orthogonality. In the case of type~5, suppose that there exists a non-trivial orthogonal decomposition of this isotropic pair into two summands, with ambient spaces $V_1$ and $V_2$ respectively. Since $\dim V = 3$, we may assume without loss of generality that $\dim V_1 = 1$. Then~$V_1$ is isotropic, and so its radical $R_1$ is equal to $V_1$. Since $R_1 = R \cap V_1$ (cf.\ Remark~\ref{radical distributivity}) and $\dim R = 1$, we f\/ind that $R = V_1$. Because~$A$,~$B$, and~$R$ are pairwise independent, it follows that $A \subseteq V_2$ and $B \subseteq V_2$ and, in turn, $A \oplus B \subseteq V_2$. But $R \subseteq A \oplus B$, which implies $R \subseteq V_2$, a contradiction to $R \subseteq V_1$.

The fact that these are all the indecomposable types follows from the decomposability to be proved next.

To prove that any isotropic pair can be decomposed as an orthogonal direct sum of indecomposables of the above types, it is enough to show that any elementary isotropic pair of type $i$,
where $i \in \{1,\dots ,10\}$,  is isomorphic to the orthogonal direct sum of copies of the $i$'th indecomposable pair above. This is straightforward when the presymplectic structure is either non-degenerate or zero (in the former case one can choose a symplectic basis adapted to the isotropic subspaces, and in the latter any choice of basis will do the job). We leave the details of this argument  to the reader and focus on the case when $(A,B)$ is of type 5.   In this case, we know that $\dim V = 3n$ for some $n \in \mathbb{N}$ and that $A \oplus B  = R \oplus A = R \oplus B$. Let $Q$
denote the subspace~$A \oplus B$ and let~$U$ be a complement of~$R$ in~$V$ such that $A \subseteq U$. For dimension reasons, $A$ is lagrangian in the symplectic subspace~$U$; let $P$ be a lagrangian complement of~$A$ in~$U$.

Because $B \subseteq R \oplus A$, if $\{b_1,\dots  ,b_n \}$ is a basis of $B$, each $b_i$ has a unique decomposition
\begin{gather}\label{b_i decomp}
b_i = r_i + a_i
\end{gather}
with  $r_i \in R \backslash \{0\}$ and $a_i \in A \backslash \{0\}$. It is a routine exercise to show that $\{r_1,\dots  ,r_n \}$ and $\{a_1,\dots  ,a_n \}$ are bases of $R$ and $A$, respectively.
Since $P$ is a lagrangian complement of $A$ in $U$, one can f\/ind a basis $\{p_1,\dots  ,p_n \}$ of $P$ which together with $\{a_1,\dots  ,a_n \}$ forms a symplectic basis of $U$.
For each $i$, set $R_i := \langle r_i \rangle$, $A_i := \langle a_i \rangle$, $B_i := \langle b_i \rangle$, $P_i := \langle p_i \rangle$, $Q_i := \langle r_i, a_i \rangle$, and $U_i := \langle a_i,p_i \rangle$.  Because $R$, $A$, $B$ and $P$ are pairwise independent, the $R_i$, $A_i$, $B_i$ and $P_i$ are too, and from (\ref{b_i decomp}) it follows that $Q_i = R_i \oplus A_i = R_i \oplus B_i = A_i \oplus B_i$. Thus, for each $i$, $(A_i,B_i)$ is an indecomposable isotropic pair of type $5$ in the presymplectic space $V_i = R_i \oplus U_i$.   From the properties of a symplectic basis, it follows that $U = U_1 \oplus \dots  \oplus U_n$ is an orthogonal decomposition, and because each $R_i$ is orthogonal to all of~$V$, the decomposition $V = V_1 \oplus \dots  \oplus V_n$ is orthogonal as well.

To show that the multiplicities of the indecomposable types in a decomposition are invariants of the pair being decomposed, we f\/irst write down a set of invariants associated to any (in general decomposable) isotropic pair $(A,B)$ in presymplectic~$V$:
\begin{alignat}{3}
&k_1 = \tfrac{1}{2}(\dim V-\dim R),\qquad  && k_6 = \dim R \cap A, &\nonumber  \\
&k_2 = \dim R,\qquad  && k_7 = \dim R \cap B, &\nonumber  \\
&k_3 =  \dim A,\qquad  && k_8 = \dim R \cap A \cap B, &  \label{canonical invariants}\\
&k_4 =  \dim B,\qquad  && k_9 = \dim R \cap (A + B), &\nonumber  \\
&k_5 =  \dim A  \cap B,\qquad  &&   k_{10} = \dim A^{\perp} \cap B. &\nonumber
\end{alignat}
Given \looseness=-1 a~decomposition of $(A,B)$ into indecomposables, one can group together all indecomposables of each type to obtain a ten-part decomposition of~$(A,B)$ into elementary types, as in Proposition~\ref{elementary decomposition}. By Lemma~\ref{distributivity properties}, this decomposition is distributive with respect to each of the subspaces in~(\ref{canonical invariants}), so each has a corresponding decomposition. Accordingly, the column vector~$\textbf{k}$ of invariants is obtained from the column vector~$\textbf{n}$ of multiplicities by multiplication by the matrix
\begin{gather*}
{\mathbf M} =
\left(
\begin{matrix}
0 & 0 & 0 & 0 & 1 & 1 & 1 & 1 & 1 & 1 \\
1 & 1 & 1 & 1 & 1 & 0 & 0 & 0 & 0 & 0 \\
0 & 1 & 1 & 0 & 1 & 1 & 1 & 0 & 1 & 0 \\
0 & 1 & 0 & 1 & 1 & 1 & 0 & 1 & 1 & 0 \\
0 & 1 & 0 & 0 & 0 & 1 & 0 & 0 & 0 & 0 \\
0 & 1 & 1 & 0 & 0 & 0 & 0 & 0 & 0 & 0 \\
0 & 1 & 0 & 1 & 0 & 0 & 0 & 0 & 0 & 0 \\
0 & 1 & 0 & 0 & 0 & 0 & 0 & 0 & 0 & 0 \\
0 & 1 & 1 & 1 & 1 & 0 & 0 & 0 & 0 & 0 \\
0 & 1 & 0 & 1 & 1 & 1 & 0 & 1 & 0 & 0
\end{matrix}
\right),
\end{gather*}
which has unit determinant and is thus invertible over $\mathbb{Z}$. Hence, the multiplicities~$\textbf{n}$ can be recovered from the invariants~$\textbf{k}$ via the inverse of $\mathbf{M}$ and are themselves invariants; i.e., decompositions of equivalent pairs into indecomposables lead to the same multiplicities.
\end{proof}

\begin{Remark}
The set of possible vectors $\textbf{n}$ of multiplicities is  of course $\mathbb{Z}_{\geq 0}^{10}$; the set of possible vectors~$\textbf{k}$ is a subset of $\mathbb{Z}_{\geq 0}^{10}$ constrained only by the ten inequalities arising from the condition that $\textbf{n} \in \mathbb{Z}_{\geq 0}^{10}$.
\end{Remark}

\begin{Remark}
The f\/irst 9 invariants in~(\ref{canonical invariants}) solve the ``triple of subspaces problem" of classifying, up to linear isomorphism, three arbitrary linear subspaces~$R$, $A$, $B$ of a vector space $V$ (without presymplectic structure). For such subspace triples, all indecomposable types are $1$-dimensional, except for one type which is $2$-dimensional. This $2$-dimensional type is closely related to the $5$th indecomposable type of isotropic pair listed in Theorem~\ref{classification}.\footnote{For the classif\/ication of subspace triples, see Brenner \cite[pp.~109--110]{Brenner_Three_Subspaces}, and Etingof et al.~\cite[pp.~84--86]{Etingof}.}

The tenth invariant $k_{10}$ in (\ref{canonical invariants}) carries the essential information  (other than $R$) connected with the presymplectic structure. For this invariant one can just as well choose $\dim A^{\perp} \cap B$; this shows, in particular, that the isotropic pairs~$(A,B)$ and~$(B,A)$ are equivalent when $\dim A = \dim B$,\footnote{This is the case, for example, in the symplectic case when the corresponding coisotropics~$A^\perp$ and~$B^\perp$ are the range and domain respectively of a linear canonical relation in $V \times \overline V.$} since all other invariants are symmetric with respect to~$A$ and~$B$.
\end{Remark}

\begin{Remark}\label{rep theory}
The elementary types listed in Proposition~\ref{elementary decomposition} are like the isotypic components in group representations.   The invariants~(\ref{canonical invariants}) are therefore, in some sense, analogous to characters, which are linear functions of the multiplicities, and which determine a representation up to equivalence in many situations. We wonder whether the present results might be recovered, for example, from the theory of some general class of quiver representations, perhaps by relations rather than maps. Note that Theorem~\ref{classification} is a Krull--Schmidt type statement, in spirit similar to Theorem~2.19 in~\cite{Etingof}, for example. It could probably also be framed in more abstract terms, such as in~\cite{Atiyah}.
\end{Remark}

The ten types of indecomposable isotropic pairs given in Theorem~\ref{classification} can be visualized as follows ($A$ is green, $B$ is blue, $R$ is red):
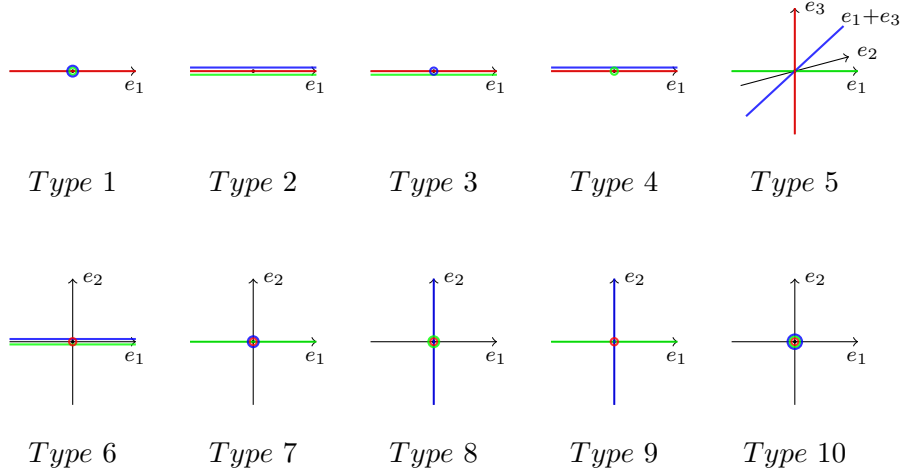
\begin{center}
\begin{tikzpicture}[scale = 2.4]
\foreach \x in {1,2,3,4,5}
	\fill[black] (\x, 2) circle (.01);
\foreach \x in {1,2,3,4,5}
	\fill[black] (\x, 0.5) circle (.01);
\foreach \x in {1,2,3,4,5}
	\draw[->] (\x - .35, 0.5) -- (\x + .35, 0.5);
\foreach \x in {1,2,3,4,5}
	\draw[->] (\x - .35, 2) -- (\x + .35, 2);
\foreach \x in {1,2,3,4,5}
	\draw[->] (\x , 0.5 - .35) -- (\x , 0.5 + .35);
\foreach \x in {5}
	\draw[->] (\x, 2  - .35) -- (\x , 2 + .35);
\foreach \x in {5}
	\draw[->] (\x - .3, 2 - .08) -- (\x + .3 , 2 + .08);
\foreach \x in {1,2,3,4,5}
	\draw (\x + .35, 0.5) node[anchor=north] {$_{e_1}$};
\foreach \x in {1,2,3,4,5}
	\draw (\x + .35, 2) node[anchor=north] {$_{e_1}$};
\foreach \x in {1,2,3,4,5}
	\draw (\x , 0.5 + .35) node[anchor=west] {$_{e_2}$};
\foreach \x in {5}
	\draw (\x + .29 , 2 + .1) node[anchor=west] {$_{e_2}$};	
\foreach \x in {5}
	\draw (\x , 2 + .35) node[anchor=west] {$_{e_3}$};
\foreach \x in {5}
	\draw (\x + .2 , 2 + .3) node[anchor=west] {$_{e_1 + e_3}$};
\foreach \x/\xtext in {1/Type \ 1,2/Type \ 2,3/Type \ 3,4/Type \ 4,5/Type \ 5}
	\draw (\x, 2 -.5) node[anchor=north] {$\xtext$};
\foreach \x/\xtext in {1/Type \ 6,2/Type \ 7,3/Type \ 8,4/Type \ 9,5/Type \ 10}
	\draw (\x , .5 - .5) node[anchor=north] {$\xtext$};
\foreach \y in {0.5}
\foreach \x in {2,4}
	\draw[green,  thick, opacity=.8] (\x - .35, \y) -- (\x + .35, \y);
\foreach \y in {0.5}
\foreach \x in {1}
	\draw[green,  thick, opacity=.8] (\x - .35, \y - 0.015) -- (\x + .35, \y - 0.015);
\foreach \y in {0.5}
\foreach \x in {1}
	\draw[blue,  thick, opacity=.8] (\x - .35, \y + 0.015) -- (\x + .35, \y + 0.015);
\foreach \y in {0.5}
\foreach \x in {3,4}
	\draw[blue,  thick, opacity=.8] (\x , \y - .35) -- (\x , \y + .35);
\foreach \y in {0.5}
\foreach \x in {1,2,3,4,5}
	\draw[red, thick, opacity=.8] (\x, \y) circle (.02);
\foreach \y in {0.5}
\foreach \x in {3,5}
	\draw[green, thick, opacity=.8] (\x, \y) circle (.03);
\foreach \y in {0.5}
\foreach \x in {2}
	\draw[blue, thick, opacity=.8] (\x, \y) circle (.03);
\foreach \y in {0.5}
\foreach \x in {5}
	\draw[blue, thick, opacity=.8] (\x, \y) circle (.04);
\foreach \y in {2}
\foreach \x in {5}
	\draw[green,  thick, opacity=.8] (\x - .35, \y) -- (\x + .35, \y);
\foreach \y in {2}
\foreach \x in {5}
	\draw[blue,  thick, opacity=.8] (\x - .27, \y  - .25) -- (\x + .27 , \y + .25);
\foreach \y in {2}
\foreach \x in {5}
	\draw[red,  thick, opacity=.8] (\x, \y  - .35) -- (\x , \y + .35);
\foreach \y in {2}
\foreach \x in {1,2,3,4}
	\draw[red,  thick, opacity=.8] (\x - .35, \y) -- (\x + .35, \y);
\foreach \y in {2}
\foreach \x in {2,3}
	\draw[green,  thick, opacity=.8] (\x - .35, \y - 0.02) -- (\x + .35, \y - 0.02);
\foreach \y in {2}
\foreach \x in {2,4}
	\draw[blue,  thick, opacity=.8] (\x - .35, \y + 0.02) -- (\x + .35, \y + 0.02);
\foreach \y in {2}
\foreach \x in {1}
	\draw[green, thick, opacity=.8] (\x, \y) circle (.02);
\foreach \y in {2}
\foreach \x in {4}
	\draw[green, thick, opacity=.8] (\x, \y) circle (.02);
\foreach \y in {2}
\foreach \x in {3}
	\draw[blue, thick, opacity=.8] (\x, \y) circle (.02);	
\foreach \y in {2}
\foreach \x in {1}
	\draw[blue, thick, opacity=.8] (\x, \y) circle (.03);
\end{tikzpicture}
\end{center}

\subsection*{Acknowledgements}

Jonathan Lorand was partially supported by ETH Zurich, the city of Zurich, and the Anna \& Hans K\"agi Foundation. Part of this research was conducted while he was at UC Berkeley as a~Visiting Student Researcher.
The authors wish to thank the referees, in particular for comments which led to a more concise presentation of our results.

\pdfbookmark[1]{References}{ref}
\LastPageEnding


\begin{thebibliography}{99}
\footnotesize \itemsep=-0.5pt

\bibitem{Atiyah}
Atiyah M.F., On the {K}rull--{S}chmidt theorem with application to sheaves,
  \textit{Bull. Soc. Math. France} \textbf{84} (1956), 307--317.

\bibitem{Benenti_Tulczyjew}
Benenti S., Tulczyjew W., Symplectic linear relations, \textit{Mem. Accad. Sci.
  Torino Cl. Sci. Fis. Mat. Natur.~(5)} \textbf{5} (1981), 71--140.

\bibitem{Brenner_Three_Subspaces}
Brenner S., Endomorphism algebras of vector spaces with distinguished sets of subspaces,
\href{http://dx.doi.org/10.1016/0021-8693(67)90016-6}{\textit{J.~Algebra}}
\textbf{6} (1967), 100--114.

\bibitem{Etingof}
Etingof P., Golberg O., Hensel S., Liu T., Schwendner A., Vaintrob D., Yudovina~E., Introduction to representation theory, {L}ecture notes, available at
  \url{http://www-math.mit.edu/~etingof/replect.pdf}.

\bibitem{Gelfand_Ponomarev}
Gel'fand I.M., Ponomarev V.A., Problems of linear algebra and classif\/ication of
  quadruples of subspaces in a f\/inite-dimensional vector space, in Hilbert
  Space Operators and Operator Algebras ({P}roc. {I}nternat. {C}onf., {T}ihany,
  1970), \textit{Colloq. Math. Soc. J\'anos Bolyai}, Vol.~5, North-Holland,
  Amsterdam, 1972, 163--237.

\bibitem{Gutt}
Gutt J., Normal forms for symplectic matrices, \href{http://dx.doi.org/10.4171/PM/1944}{\textit{Port. Math.}} \textbf{71}
  (2014), 109--139, \href{http://arxiv.org/abs/1307.2403}{arXiv:1307.2403}.

\bibitem{Li-Bland_Weinstein}
Li-Bland D., Weinstein A., Selective categories and linear canonical relations,
  \href{http://dx.doi.org/10.3842/SIGMA.2014.100}{\textit{SIGMA}} \textbf{10} (2014), 100, 31~pages, \href{http://arxiv.org/abs/1401.7302}{arXiv:1401.7302}.

\bibitem{Lorand}
Lorand J., Classifying linear canonical relations, \href{http://arxiv.org/abs/1508.04568}{arXiv:1508.04568}.

\bibitem{Lorand_Weinstein}
Lorand J., Weinstein A., Coisotropic pairs, \href{http://arxiv.org/abs/1408.5620}{arXiv:1408.5620}.

\bibitem{lo-we:decomposition}
Lorand J., Weinstein A., Decomposition of (co)isotropic relations, {i}n preparation.

\bibitem{Roman}
Roman S., Advanced linear algebra, \href{http://dx.doi.org/10.1007/978-0-387-72831-5}{\textit{Graduate Texts in Mathematics}}, Vol.~135, 3rd ed., Springer, New York, 2008.

\bibitem{Sergeichuk}
Sergeichuk V.V., Classif\/ication of pairs of subspaces in spaces with scalar
  product, \href{http://dx.doi.org/10.1007/BF01071340}{\textit{Ukrainian Math.~J.}} \textbf{42} (1990), 487--491.

\bibitem{Towber}
Towber J., Linear relations, \href{http://dx.doi.org/10.1016/0021-8693(71)90112-8}{\textit{J.~Algebra}} \textbf{19} (1971), 1--20.

\bibitem{we:coisotropic}
Weinstein A., Coisotropic calculus and {P}oisson groupoids, \href{http://dx.doi.org/10.2969/jmsj/04040705}{\textit{J.~Math.
  Soc. Japan}} \textbf{40} (1988), 705--727.

\bibitem{Weinstein}
Weinstein A., A note on the {W}ehrheim--{W}oodward category, \href{http://dx.doi.org/10.3934/jgm.2011.3.507}{\textit{J.~Geom.
  Mech.}} \textbf{3} (2011), 507--515, \href{http://arxiv.org/abs/1012.0105}{arXiv:1012.0105}.

\bibitem{Weinstein2}
Weinstein A., Categories of (co)isotropic linear relations, \href{http://arxiv.org/abs/1503.06240}{arXiv:1503.06240}.

\bibitem{Williamson}
Williamson J., On the normal forms of linear canonical transformations in
  dynamics, \href{http://dx.doi.org/10.2307/2371583}{\textit{Amer.~J. Math.}} \textbf{59} (1937), 599--617.

\end{thebibliography}
\end{document}